\documentclass[12pt]{article}
\usepackage{amssymb,amsmath,,amsthm,tikz}
\newtheorem{theorem}{Theorem}

\newtheorem{question}[theorem]{Quesion}

\title{Matching integral 2-connected graphs}

\author{ Ebrahim Ghorbani
\\[.3cm]
{\sl\normalsize Department of Mathematics, K. N. Toosi University of Technology,}\\
{\sl\normalsize P. O. Box 16765-3381, Tehran, Iran}}

\begin{document}
\maketitle
\footnotetext{{\em E-mail Address}:  {\tt e\_ghorbani@ipm.ir} }

\begin{abstract}
In this note we present an infinite family of 2-connected graphs such that their matching polynomials have only integer zeros.
This answers in negative a question of Akbari et al. [Graphs with integer matching polynomial zeros, Discrete Appl. Math. 224 (2017), 1--8].

\vspace{5mm}
\noindent {\bf Keywords:} Matching polynomial, Matching integral, 2-connected graph  \\[.1cm]
\noindent {\bf AMS Mathematics Subject Classification\,(2010):} 05C31, 05C40
\end{abstract}

Let $G$ be a graph.
An {\em $r$-matching}  in $G$
is a set of $r$ pairwise non-incident edges. The number of $r$-matchings of $G$ is denoted by $m(G, r)$. The  {\em matching polynomial} of $G$ is defined by $$\mu(G, x)=\sum_{r=0}^{\lfloor{\frac{n}{2}}\rfloor}(-1)^{r}m(G, r)x^{n-2r},$$ where $n$ is the order of $G$ and $m(G, 0)$ is considered to be $1$. For basic properties of matching polynomials see \cite{gg,gu}.

A graph is called {\em matching integral} if the zeros of its matching polynomial are all integers.
In \cite{acgk}, matching integral graphs within some specific families of graphs including connected regular graphs, claw-free graphs and graphs with perfect matching were characterized.  In particular, the following question was asked:

\begin{question}[Akbari et al. \cite{acgk}]\rm
Is it true that there are finitely many matching integral 2-connected graphs?
\end{question}
 We answer this question in negative:
\begin{theorem}
There are infinitely many matching integral $2$-connected graphs.
\end{theorem}

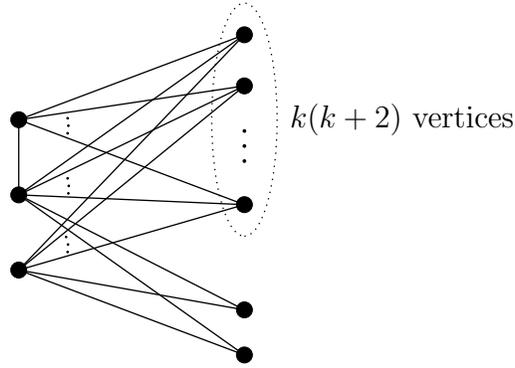
\begin{figure}
\centering\begin{tikzpicture}
\draw [line width=.5pt] (0,2)-- (0,1);
\draw [line width=.5pt] (3,-1.13)-- (0,0);
\draw [line width=.5pt] (3,-1.13)-- (0,1);
\draw [line width=.5pt] (3,-0.53)-- (0,0);
\draw [line width=.5pt] (3,-0.53)-- (0,1);
\draw [line width=.5pt] (3,0.87)-- (0,0);
\draw [line width=.5pt] (3,0.87)-- (0,1);
\draw [line width=.5pt] (3,0.87)-- (0,2);
\draw [line width=.5pt] (3,3.13)-- (0,2);
\draw [line width=.5pt] (3,3.13)-- (0,1);
\draw [line width=.5pt] (3,3.13)-- (0,0);
\draw [line width=.5pt] (3,2.45)-- (0,2);
\draw [line width=.5pt] (3,2.45)-- (0,1);
\draw [line width=.5pt] (3,2.45)-- (0,0);
\draw [dotted,rotate around={-90:(3,2)},line width=.5pt] (3,2) ellipse (1.55cm and 0.435cm);
\draw [fill=black] (3,1.85) circle (.5pt);
\draw [fill=black] (3,1.65) circle (.5pt);
\draw [fill=black] (3,1.45) circle (.5pt);
\draw (5.1,2) node {$k(k+2)$ vertices};
\draw [fill=black] (.65,2.02) circle (.3pt);
\draw [fill=black] (.66,1.92) circle (.3pt);
\draw [fill=black] (.65,1.82) circle (.3pt);
\draw [fill=black] (.65,1.22) circle (.3pt);
\draw [fill=black] (.665,1.12) circle (.3pt);
\draw [fill=black] (.66,1.02) circle (.3pt);
\draw [fill=black] (.63,.44) circle (.3pt);
\draw [fill=black] (.655,.34) circle (.3pt);
\draw [fill=black] (.65,.24) circle (.3pt);
\draw [fill=black] (0,0) circle (3pt);
\draw [fill=black] (0,2) circle (3pt);
\draw [fill=black] (0,1) circle (3pt);
\draw [fill=black] (3,-0.53) circle (3pt);
\draw [fill=black] (3.,-1.13) circle (3pt);
\draw [fill=black] (3,0.87) circle (3pt);
\draw [fill=black] (3,3.13) circle (3pt);
\draw [fill=black] (3,2.45) circle (3pt);
\end{tikzpicture}
\caption{The graph $H_k$}\label{fig:Hk}
\end{figure}

\begin{proof}{For any positive integer $k$, we define the graph $H_k$ as depicted in Figure~\ref{fig:Hk}.
Clearly, $H_k$ is $2$-connected.
It can be easily verified that
\begin{align*}
    m(H_k,1)&=3k(k+2)+5\\
    &=k^2+(k+1)^2+(k+2)^2,\\
    m(H_k,2)&=\big(k(k+2)+1\big)\big(3k(k+2)+2\big)+k(k+2)+2\\
    &=k^2(k+1)^2+k^2(k+2)^2+(k+1)^2(k+2)^2,\\
    m(H_k,3)&=\big(k(k+2)\big)^2\big(k(k+2)+1\big)\\
    &=k^2(k+1)^2(k+2)^2,
\end{align*}
and $m(H_k,4)=0$. It follows that
\begin{align*}
\mu(H_k,x)&=x^{k(k+2)-1}\left(x^6-m(H_k,1)x^4+m(H_k,2)x^2-m(H_k,3)\right)\\
&=x^{k(k+2)-1}\left(x^2-k^2\right)\left(x^2-(k+1)^2\right)\left(x^2-(k+2)^2\right).
\end{align*}
This shows that $H_k$ is a matching integral graph.
}\end{proof}

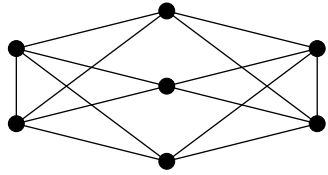
\begin{figure}
\centering\begin{tikzpicture}
\draw [line width=.5pt] (-2,0.5)-- (-2,-0.5);
\draw [line width=.5pt] (0,1)-- (-2,0.5);
\draw [line width=.5pt] (0,1)-- (-2,-0.5);
\draw [line width=.5pt] (0,0)-- (-2,0.5);
\draw [line width=.5pt] (0,0)-- (-2,-0.5);
\draw [line width=.5pt] (0,-1)-- (-2,-0.5);
\draw [line width=.5pt] (0,-1)-- (-2,0.5);
\draw [line width=.5pt] (0,1)-- (2,0.5);
\draw [line width=.5pt] (0,0)-- (2,0.5);
\draw [line width=.5pt] (0,-1)-- (2,0.5);
\draw [line width=.5pt] (0,1)-- (2,-0.5);
\draw [line width=.5pt] (0,0)-- (2,-0.5);
\draw [line width=.5pt] (0,-1)-- (2,-0.5);
\draw [line width=.5pt] (2,0.5)-- (2,-0.5);
\draw [fill=black] (0,0) circle (3pt);
\draw [fill=black] (0,1) circle (3pt);
\draw [fill=black] (0,-1) circle (3pt);
\draw [fill=black] (-2,0.5) circle (3pt);
\draw [fill=black] (-2,-0.5) circle (3pt);
\draw [fill=black] (2,0.5) circle (3pt);
\draw [fill=black] (2,-0.5) circle (3pt);
\end{tikzpicture}
\caption{A 3-connected  matching integral graph}\label{fig:3-con}
\end{figure}

A natural generalization would be looking for matching integral graphs with higher connectivity.
We are only known of one $3$-connected matching integral graph, namely the graph depicted in Figure~\ref{fig:3-con}, whose matching polynomial is
$$x^7-14x^5+49x^3-36x=x(x^2-1)(x^2-4)(x^2-9).$$

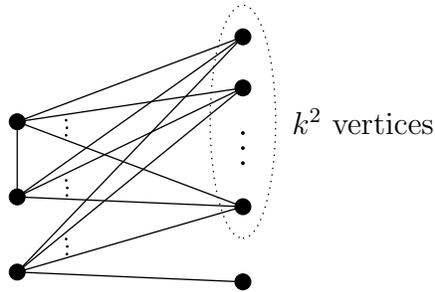
\begin{figure}[h!]
\centering\begin{tikzpicture}
\draw [line width=.5pt] (0,2)-- (0,1);
\draw [line width=.5pt] (3,-.13)-- (0,0);
\draw [line width=.5pt] (3,0.87)-- (0,0);
\draw [line width=.5pt] (3,0.87)-- (0,1);
\draw [line width=.5pt] (3,0.87)-- (0,2);
\draw [line width=.5pt] (3,3.13)-- (0,2);
\draw [line width=.5pt] (3,3.13)-- (0,1);
\draw [line width=.5pt] (3,3.13)-- (0,0);
\draw [line width=.5pt] (3,2.45)-- (0,2);
\draw [line width=.5pt] (3,2.45)-- (0,1);
\draw [line width=.5pt] (3,2.45)-- (0,0);
\draw [dotted,rotate around={-90:(3,2)},line width=.5pt] (3,2) ellipse (1.55cm and 0.435cm);
\draw [fill=black] (3,1.85) circle (.5pt);
\draw [fill=black] (3,1.65) circle (.5pt);
\draw [fill=black] (3,1.45) circle (.5pt);
\draw (4.6,2) node {$k^2$ vertices};
\draw [fill=black] (.65,2.02) circle (.3pt);
\draw [fill=black] (.66,1.92) circle (.3pt);
\draw [fill=black] (.65,1.82) circle (.3pt);
\draw [fill=black] (.65,1.22) circle (.3pt);
\draw [fill=black] (.665,1.12) circle (.3pt);
\draw [fill=black] (.66,1.02) circle (.3pt);
\draw [fill=black] (.65,.44) circle (.3pt);
\draw [fill=black] (.665,.34) circle (.3pt);
\draw [fill=black] (.66,.24) circle (.3pt);
\draw [fill=black] (0,0) circle (3pt);
\draw [fill=black] (0,2) circle (3pt);
\draw [fill=black] (0,1) circle (3pt);
\draw [fill=black] (3,-.13) circle (3pt);
\draw [fill=black] (3,0.87) circle (3pt);
\draw [fill=black] (3,3.13) circle (3pt);
\draw [fill=black] (3,2.45) circle (3pt);
\end{tikzpicture}
\caption{The graph $H'_k$}\label{fig:H'k}
\end{figure}

Our search to find more $3$-connected matching integral graphs was not successful. However,
we may present another family of  matching integral graphs, namely the graphs $H'_k$ (see Figure~\ref{fig:H'k}).
It is straightforward to check that
$$\mu(H'_k,x)=x^{k^2-2}\left(x^2-(k-1)^2\right)\left(x^2-k^2\right)\left(x^2-(k+1)^2\right).$$
 $H'_k$ has a cut vertex but for $k\ge2$ it is very close to being $3$-connected in the sense that it can be obtained
by attaching a pendant vertex to a $3$-connected graphs.

These observations motivate us to put forward the following question.

\begin{question}\rm Are there matching integral $k$-connected graphs for integers $k\ge4$?
\end{question}

\section*{Acknowledgments}
This work was done while the author was visiting Alfred R\'enyi Institute  of Mathematics in August 2017. He would like to thank
the institute for its hospitality and support and Peter Csikv\'ari for fruitful conversations.

\end{document}